\documentclass{amsart}
\usepackage{amsmath}
\usepackage{amsmath}
\usepackage{amssymb}

\input xy
\xyoption{all}

\newcommand{\free}[1]{\,{\displaystyle{\ast}}_{#1}\,}

\title[Exactness of finite dimensional $C^*$-algebras]{Exactness of universal
free products of finite dimensional $C^*$-algebras with
amalgamation}

\author{Benton L. Duncan}

\address{Department of Mathematics\\
North Dakota State University\\
Fargo, ND  \\
USA}

\email{benton.duncan@ndsu.edu}



\begin{document}

\theoremstyle{plain}
\newtheorem{theorem}{Theorem}
\newtheorem{lemma}{Lemma}
\newtheorem{proposition}{Proposition}
\newtheorem{corollary}{Corollary}

\theoremstyle{definition}
\newtheorem{dfn}{Definition}
\newtheorem*{construction}{Construction}
\newtheorem{example}{Example}

\theoremstyle{remark}
\newtheorem*{conjecture}{Conjecture}
\newtheorem*{acknowledgement}{Acknowledgements}
\newtheorem{rmk}{Remark}

\setcounter{MaxMatrixCols}{16}

\begin{abstract} We investigate free products of finite dimensional
$C^*$-algebras with amalgamation over diagonal subalgebras.  We look
to determine under what circumstances a given free product is exact
and/or nuclear.  In some cases we find a description of the algebra
in terms of a more readily understood algebra.
\end{abstract}

\maketitle

Recall that for $C^*$-algebras $A$ and $B$ there are many possible
norms on $A \otimes B$ for which the completion is a $C^*$-algebra.
In particular there are two standard completions $A \otimes_{\rm
min} B$ and $A \otimes_{\rm max} B$ corresponding to the `smallest'
and `largest' possible tensor norms. We say that $A$ is nuclear if
these two tensor products correspond for all $C^*$-algebras $B$. We
say that a $C^*$-algebra $A$ is exact if given any short exact
sequence of $C^*$-algebras \[ 0 \rightarrow C \rightarrow B
\rightarrow B/C \rightarrow 0 \] the associated sequence \[ 0
\rightarrow C \otimes_{\rm min} A \rightarrow B \otimes_{\rm min} A
\rightarrow B/C \otimes_{\rm min} A \rightarrow 0 \] is a short
exact sequence.  Both of these properties represent some appreciable
level of `niceness' for $C^*$-algebras. (For more information about
nuclear and exact $C^*$-algebras we refer the reader to
\cite{Blackadar}).

In this paper we are interested in the question of whether the
universal free product of matrix algebras, with and without
amalgamation over diagonal subalgebras, are exact and/or nuclear.
This question is solved in the case of the reduced free product, see
\cite{Dykema} where it is shown that the reduced amalgamted free
product of exact $C^*$-algebras is exact.  Of course the universal
free products are often not `nice' in any reasonable sense; this is
borne out in this paper by the fact that even simple finite
dimensional $C^*$-algebras (matrix algebras over $\mathbb{C}$)
quickly lose exactness and/or nuclearity when dealing with free
products. However there are cases in which nuclearity is preserved
under universal free products. This work was motivated by
\cite{Duncan} where the question of nuclearity/exactness was
discussed for free products of directed graph $C^*$-algebras.  Since
directed graph $C^*$-algebras are often free products of finite
dimensional $C^*$-algebras this paper was the natural outgrowth of
that investigation.

This work is related although different from \cite{Albeverio,
Jushenko} where a related notion of $*$-wildness for finite
dimensional free products was discussed.  There is an important
distinction between the present investigation and the aforementioned
work: the free products in \cite{Albeverio, Jushenko} are all
assumed to be unital.  We look at a broader class of possible free
products, allowing amalgamations over different diagonal
subalgebras.

Unless specifically stated otherwise an algebra in this paper will
be a $C^*$-algebra, and an isomorphism of algebras will mean a
$*$-isomorphism.  By $M_j$ we mean the $j\times j$ matrices over
$\mathbb{C}$; for the purposes of this paper we will always assume
that $1 < j < \infty$. By a unital diagonal subalgebra of $M_j$ we
mean a subalgebra of the $j \times j$ diagonal matrices which
contains the unit of $M_j$.  The notation $A \free{}B$ will denote
the universal free product of the algebras $A$ and $B$ with no
amalgamation.  When $A$ and $B$ contain a common subalgebra $D$ we
will write $A \free{D} B$ to denote the universal free product of
$A$ and $B$ with amalgamation over $D$.

\section{Amalgamation diagrams}

Given $M_j$ and $M_k$ we intend to look at algebras of the form $M_j
\free{D} M_k$ where $D$ is a copy of $\mathbb{C}^n$ embedded into
the two algebras as a diagonal subalgebra of the matrix algebras. Of
course there are many different ways to do this embedding.  We
introduce some notation to describe how $\mathbb{C}^n$ embeds into
$M_j$.

We use the diagram \begin{align*} & M_j: \boxed{j_1} \boxed{j_2}
\boxed{j_3} \cdots \boxed{j_n} \boxed{0}
\end{align*} to describe the embedding \[ \begin{bmatrix} \lambda_1 \\ \lambda_2
\\ \vdots \\ \lambda_n \end{bmatrix} \mapsto  \begin{bmatrix} \lambda_1 I_{j_1} & 0 &
0 & \cdots & 0 \\ 0 & \lambda_2I_{j_2} & 0 & \cdots & 0 \\
\vdots & \vdots & \ddots & \vdots & \vdots \\ 0 & 0 & \cdots &
\lambda_{j_n}I_{j_n} & 0 \\ 0 & 0 & \cdots & 0 & 0
\end{bmatrix}  \in \begin{bmatrix} M_{j_1} & 0 & 0 & \cdots & 0 \\ 0 & M_{j_2} & 0
& \cdots & 0 \\ \vdots & \vdots & \ddots & \vdots & \vdots
\\ 0 & \cdots & 0 & M_{j_n} & 0 \\ 0 & \cdots & 0 & 0 & 0
\end{bmatrix} \subseteq M_j. \] Here $I_{\alpha}$ is the $\alpha \times \alpha$
identity matrix in $M_{\alpha}$. If there is a box containing a zero
then we call such a box a {\em zero-box}.  Further, in our notation
there is at most one zero box.

Notice that the embedding diagram tells us:
\begin{enumerate} \item the value of $n$,
\item and whether the embedding is non-unital, indicated by the presence
of a zero-box. \end{enumerate}  For the purposes of our results it
is safe to assume that through the use of elementary row operations
that any zero-box is listed last.

Now when looking at the free product of two matrix algebras with
amalgamation over $\mathbb{C}^n$ it is clear that just writing $M_j
\free{\mathbb{C}^n} M_k$ will be unsuitable because it is not clear
how we are embedding $\mathbb{C}^n$ into the two matrix algebras. To
see the amalgamation we will use pairs of embedding diagrams.  We
will call a pair of embedding diagrams an {\em amalgamation diagram}
since they represent the amalgamating subalgebra in a free product.
We will present two examples to illustrate how this will work.

\begin{example}\label{m2on} We start with an example from W.\ Paschke
\cite[Example 3.3]{Brownexact}. There it is noted that with suitable
amalgamation $M_{j+1} \free{\mathbb{C}^2} M_2$ is isomorphic to
$M_{j+1} \otimes \mathcal{O}_j$.  The amalgamation can be described
using the
amalgamation diagram \begin{align*} & M_j:\boxed{1} \boxed{j-1} \\
 & M_2:\boxed{1} \boxed{1}. \end{align*}  Here the $M_j$-row represents
$\mathbb{C}^2$ as the subalgebra of $M_j$ given by
\[ \begin{bmatrix} \lambda_1 & 0 & \cdots & 0 \\ 0 & \lambda_2 & \cdots & 0
\\ \vdots & \vdots & \ddots & \vdots \\ 0 & \cdots & 0 &
\lambda_2 \end{bmatrix}. \]  The $M_2$-row represents the usual
embedding of $\mathbb{C}^2$ as the diagonal subalgebra of $M_2$.
Notice that for both embeddings $\mathbb{C}^2$ is a unital
subalgebra of the associated algebra.
\end{example}

\begin{example}\label{mntensora} The next example is from
\cite[Chapter 6]{Loring}. There it is shown that for unital $A$ and
appropriate choice of embedding we have that $ M_j \free{\mathbb{C}}
A $ is isomorphic to $M_j(A)$. For our notation we will look at the
specific case of $A = M_k$ and
amalgamation diagram \begin{align*} & M_j:\boxed{1} \boxed{0} \\
& M_k:\boxed{k} . \end{align*}  Here the scalar multiples of the
identity in $M_k$ are matched up with the $1\times 1$ entry in
$M_j$. \end{example}

\begin{example}\label{Mk} Finally we have the following example which,
while not of the form described above will allow us to make some
computations later.  The algebra \[ A:= \begin{bmatrix} M_k & 0 \\
0 & \mathbb{C} \end{bmatrix} \free{\mathbb{C}^{k+1}} \begin{bmatrix}
\mathbb{C}^{k-1} & 0 \\ 0 & M_2 \end{bmatrix}  \] is isomorphic to $
M_{k+1}$, where $\mathbb{C}^{k+1}$ is the canonical inclusion as
diagonal matrices. Certainly there is an onto $*$-representation
$\pi: A \rightarrow
M_{k+1}$ induced by the inclusions \begin{align*} \iota_{k,1}: \begin{bmatrix} M_k & 0 \\
0 & \mathbb{C} \end{bmatrix} & \subseteq M_{k+1} \\ \iota_{k-1,2}
\begin{bmatrix} \mathbb{C}^{k-1} & 0 \\ 0 & M_2 \end{bmatrix}
& \subseteq M_{k+1}. \end{align*}  We need only show that $M_{k+1}$
satisfies the requisite universal property.  So let $B$ be a
$C^*$-algebra and assume that we have $*$-representations $\pi_1:
\begin{bmatrix} M_k & 0 \\ 0 & \mathbb{C} \end{bmatrix}
\rightarrow B$ and $ \pi_2: \begin{bmatrix} \mathbb{C}^{k-1} & 0 \\
0 & M_2 \end{bmatrix} \rightarrow B$ with $\pi_1|_D = \pi_2|_D$ for
the subalgebra of diagonal matrices $D$.  Then for the elementary
matrices $e_{i,j} \in M_{k+1}$ define
\[\pi(e_{i,j}) = \begin{cases} \pi_1(e_{i,j}) & 1 \leq i, j < k+1 \\
\pi_2(e_{k+1,k+1}) = \pi_1(e_{k+1,k+1}) & i = j = k+1 \\
\pi_1(e_{i,k})\pi_2(e_{k,k+1}) & 1 \leq i < k+1, j = k+1 \\
\pi_2(e_{k+1,k})\pi_1(e_{k,j} & 1 \leq j < k+1, i = k+1
\end{cases}.\]  Notice that in the second case since $\pi_1|_D =
\pi_2|_D$ which tells us that the second case is well defined.  For
general matrices we extend using linearity. We need only show that
$\pi$ induces a $*$-representation on $M_{k+1}$. To verify this we
notice first that $\pi$ is linear by construction. Next, to show
that $\pi(A^*) = \pi(A)^*$ we only need show, by linearity, that $
\pi({e_{i,j}}^*) = \pi(e_{i,j})^*$ for all $i, j$. This is trivial
if $1 \leq i,j \leq k$ or $ i=j=k+1$ since $\pi_1$ and $\pi_2$ are
$*$-representations. So assume that $i < k+1$ and consider
\begin{align*} \pi(e_{i,k+1})^* & = (\pi_1(e_{i,k})
\pi_2(e_{k,k+1}))^* \\ & = \pi_2(e_{k,k+1})^*\pi_1(e_{i,k+1})^* \\
&= \pi_2(e_{k+1,k})\pi_1(e_{k+1,i}) \\ &= \pi(e_{k+1,i}) =
\pi({e_{i,k+1}}^*). \end{align*} The third equality follows since
$\pi_1$ and $\pi_2$ are $*$-representations.  The case where $j <
k+1$ and $i=k+1$ is similar.

We next need to show that $\pi$ is multiplicative.  We will consider
products of the form $e_{i,m} = e_{i,j}e_{j,m}$.  Again this follows
using cases.  If $1 \leq i,j,m < k+1$, or $i=j=m=k+1$ then
$\pi(e_{i,m}) = \pi_1(e_{i,m})= \pi_1(e_{i,j} e_{j,m}) =
\pi_1(e_{i,j}) \pi_1(e_{j,m}) = \pi(e_{i,j})\pi(e_{j,m})$. There are
six remaining cases, we will do two of them, the remainder will
follow in a similar fashion.

Assume that $ j=k+1$ and $1 \leq i, m<k+1$, then
\begin{align*} \pi(e_{i,j}) \pi(e_{j,m}) & = \pi_1(e_{i,k}) \pi_2(e_{k,k+1})
\pi_2(e_{k+1,k})\pi_1(e_{k,m}) \\ &= \pi_1(e_{i,k})
\pi_2(e_{k,k})\pi_1(e_{k,m}) \\ &= \pi_1(e_{i,k})
\pi_1(e_{k,k})\pi_1(e_{k,m}) \\ &= \pi_1 (e_{i,k}e_{k,k}e_{k,m}) =
\pi_1(e_{i,m}) = \pi(e_{i,m}). \end{align*} Notice that in the third
equality we used that $\pi_2$ is a homomorphism, in the next line we
used that $\pi_1|_D = \pi_2|_D$, and then in the line after we use
the fact that $\pi_1$ is a homomorphism.

Next consider the case that $i = k+1, m=k+1$ and $1 \leq j < k+1$
and compute \begin{align*} \pi(e_{i,j})\pi(e_{j,m}) &=
\pi_2(e_{k+1,k}) \pi_1(e_{k,j}) \pi_1(e_{j,k}) \pi_2(e_{k,k+1}) \\
&= \pi_2(e_{k+1,k} \pi_1(e_{k,k}) \pi_2(e_{k,k+1}) \\ &=
\pi_2(e_{k+1,k} \pi_2(e_{k,k}) \pi_2(e_{k,k+1}) \\ &=
\pi_2(e_{k+1,k}e_{k,k}e_{k,k+1})\\ & = \pi_2(e_{k+1,k+1}) =
\pi(e_{i,j}e_{j,m}). \end{align*}

Similar calculations finish the remaining cases and then applying
linearity completes the proof that $M_{k+1}$ has the requisite
universal property and hence is isomorphic to $A$.
\end{example}

The following will be useful in analyzing exactness and nuclearity
for free products.

\begin{theorem}\label{inclusions} If $D$ is a $C^*$-subalgebra of $A_1$ and $A_2$ then
there exists a canonical onto $*$-representation $ \pi: A_1
\free{}A_2 \rightarrow A_1 \free{D}A_2$.  If, in addition, $C$ is a
$C^*$-algebra with $D \subseteq C \subseteq A_i$ for each $i = 1,2$
then there is a canonical onto $*$-representation $\sigma: A_1
\free{D}A_2 \rightarrow A_1 \free{C} A_2$.
\end{theorem}

\begin{proof} Let $\iota_i: A_i \rightarrow A_1\free{D}A_2 $ be the
canonical inclusion (i.e. $A_i \subseteq A_1 \free{D}A_2$). Then by
the universal property of $A_1 \free{}A_2$ there exists a canonical
$*$-representation $\iota_1 \free{} \iota_2: A_1 \free{}A_2
\rightarrow A_1 \free{D} A_2$.  This map is onto since a generating
set for $A_1 \free{D} A_2$ is contained in the image of $\iota_1
\free{} \iota_2$.

Next let $\beta_i: A_i \rightarrow A_1 \free{C} A_2$ be the
canonical inclusion.  Notice that $ \beta_1(d) = \beta_2(d)$ for all
$d \in D$ since $D \subseteq C$ and hence there is an induced
$*$-representation $\beta_1 \free{D} \beta_2: A_1 \free{D} A_2
\rightarrow A_1 \free{C} A_2$ which is onto for the same reason as
the previous map. \end{proof}

The following is immediate since both nuclearity and exactness pass
to quotients.

\begin{corollary}\label{subalgebra} If $A \free{} B$ is nuclear so is $A \free{D} B$
for any $C^*$-algebra $D$ with $D \subseteq A$ and $D \subseteq B$.
If $A \free{D} B$ is not exact then neither $A \free{} B$ nor $A
\free{C} B$ is exact for any subalgebra $C$ with $D \subseteq C
\subseteq A$ and $ D \subseteq C \subseteq B$.
\end{corollary}

Finally we have one more well-known example which will provide a
standard non-exact $C^*$-algebra for our results.

\begin{example}\label{ctct} The algebra $C(\mathbb{T}) \free{\mathbb{C}}
C(\mathbb{T})$ is isomorphic to the non-exact $C^*$-algebra
$C^*(\mathbb{Z}) \free{\mathbb{C}} C^*(\mathbb{Z}) = C^*(\mathbb{Z}
\free{} \mathbb{Z}) = C^*(F_2)$ (see \cite{Wassermann} for a proof
that the latter is not exact).
\end{example}

\section{Algebras of the form $M_j \free{D} M_k$}

We have already seen two examples of these type of algebras, both of
which were nuclear.  The general case will be more complicated and
will depend on the nature of $D$, and on the embedding diagrams for
$D \subseteq M_j$ and $D \subseteq M_k$.

\begin{proposition} The algebra $M_3 \free{\mathbb{C}^3} M_3$ is
isomorphic to $M_3 \otimes (C(\mathbb{T}) \free{\mathbb{C}}
C(\mathbb{T}))$ and hence is not exact. \end{proposition}

\begin{proof} We know from Example \ref{Mk} that $M_3 = (M_2 \oplus
\mathbb{C}) \free{\mathbb{C}^3} (\mathbb{C} \oplus M_2)$ and hence
$M_3 \free{\mathbb{C}^3} M_3$ can be rewritten as \[ \left(( M_2
\oplus \mathbb{C}) \free{\mathbb{C}^3} (\mathbb{C} \oplus M_2
 )\right) \free{\mathbb{C}^3} \left( (M_2 \oplus \mathbb{C})
\free{\mathbb{C}^3} (\mathbb{C} \oplus M_2) \right). \]  Of course
by rearranging we can rewrite this as \[ \left( (M_2 \oplus
\mathbb{C}) \free{\mathbb{C}^3} (M_2 \oplus \mathbb{C}) \right)
\free{\mathbb{C}^3} \left( (\mathbb{C} \oplus M_2)
\free{\mathbb{C}^3} (\mathbb{C} \oplus M_2) \right)\] which by
Example \ref{m2on} is isomorphic to $\left(M_2(C(\mathbb{T})) \oplus
\mathbb{C}\right) \free{\mathbb{C}^3} \left(\mathbb{C} \oplus
M_2(C(\mathbb{T})) \right)$.  The latter algebra has a canonical
representation onto a generating set for the algebra $ M_3 \otimes
(C(\mathbb{T}) \free{\mathbb{C}} C(\mathbb{T}))$ via the inclusion
maps. It is a simple matter to see that the algebra $ M_3 \otimes
(C(\mathbb{T}) \free{\mathbb{C}} C(\mathbb{T}))$ satisfies the
universal property for \[\left(M_2(C(\mathbb{T})) \oplus \mathbb{C}
\right) \free{\mathbb{C}^3} \left(\mathbb{C} \oplus
M_2(C(\mathbb{T}))\right). \] Lack of exactness now follows using
Example \ref{ctct} \end{proof}

\begin{proposition}\label{mjmk} If $D$ is a unital diagonal subalgebra of $M_j$
and $M_k$ such that $\dim{D} \geq 3$, then $M_j \free{D} M_k$ is not
exact. \end{proposition}

\begin{proof} First notice that there is a $3$-dimensional
subalgebra of $D$, call it $E$ and denote by $e$ the identity of
$E$. Then set $A = \{x \in M_j: xe = ex = x \}$ and $ B = \{ x \in
M_k: xe=ex=x \}$. It is routine to verify that $ A \cong M_3$ and $B
\cong M_3$. Then applying \cite[Proposition
2.4]{Armstrong-Dykema-Exel-Li:2003} with the canonical conditional
expectations given by projections onto the appropriate subalgebras
$M_3$ we have that $ A \free{E} B \subset M_j \free{D} M_k$ and
hence $M_j \free{D} M_k$ is not exact.
\end{proof}

We let $m_i$ denote the minimum value in the $i$th row of the
amalgamation diagram for $M_j \free{\mathbb{C}^k} M_k$.  Define the
{\em minimum value of the diagram} to be the sum of the $m_i$ as $i$
ranges over each row of the amalgamation diagram.

\begin{proposition}\label{submjmk} Let $D$ be a unital diagonal subalgebra of $M_j$
and $M_k$. If the minimum value of the amalgamation diagram for $M_j
\free{D} M_k$ is greater than or equal to $3$, then $M_j \free{D}
M_k$ is not exact.
\end{proposition}

\begin{proof} By hypothesis there exist diagonal subalgebras
$E \subseteq M_j$ and $E \subseteq M_k$ such that the dimension of
$E$ is greater than or equal to $3$.  As in the previous proposition
we let $e$ denote the identity in $E$ and set $ A = \{x \in M_j: xe
= ex = x \}$ and set $B = \{x \in M_k: xe = ex = x \}$.  Notice that
$ A \cong M_{\dim E} \cong B$. The result will follow by applying
\cite[Proposition 2.4]{Armstrong-Dykema-Exel-Li:2003} with the
canonical conditional expectations and noting that $ A \free{E} B
\subseteq M_j \free{E} M_k$.\end{proof}

\begin{theorem}\label{graph} Let $D$ be a unital diagonal subalgebra of $M_j$ and
$M_k$ such that the minimum value of the amalgamation diagram for
$M_j \free{D} M_k$ is $2$.  If $\dim D = 2$ then the algebra is
nuclear. \end{theorem}

\begin{proof} We will show that such an algebra is a directed graph
$C^*$-algebra and hence is nuclear. Let $G$ be the directed graph
with $2$-vertices $\{ v_1, v_2 \}$ and $(j-1)+(k-1)$ edges $\{ e_1,
e_2, \cdots, e_{j-1}, f_1, f_2, \cdots, f_{k-1} \}$ with $r(e_i) =
v_1$, $s(e_i) = v_2$ and $ r(f_i) = v_2, s(f_i) = v_1$. We claim
that $C^*(G)$ is isomorphic to $M_j \free{D} M_k$. Notice that
$e_{i+1,1} \in M_j$ and $e_{j,k} \in M_k$ form a collection of
partial isometries which form a Cuntz-Krieger family for the graph
$G$. Further notice that this Cuntz-Krieger family generates the
algebra $M_j \free{D} M_k$.  By \cite[Proposition
1.21]{Raeburn:2006} there is a $*$-representation $\pi: C^*(G)
\rightarrow M_j \free{D} M_k$.  Notice further that the directed
graph thus constructed is cofinal and every cycle has an entry hence
$C^*(G)$ is simple by \cite[Proposition 4.2]{Raeburn:2006}.  It
follows that $\pi$ is one-to-one and hence the free product algebra
is isomorphic to $C^*(G)$ and is nuclear.\end{proof}

Notice that the minimum value of an amalgamation diagram for $M_j
\free{D} M_k$ is never equal to $1$.  For unital amalgamations of
finite dimensional algebras we have one case remaining.

\begin{proposition} The algebra $M_2 \free{\mathbb{C}} M_2$ is not exact.
\end{proposition}

\begin{proof} Define $\pi_1: M_2 \rightarrow M_2(C(\mathbb{T})
\free{\mathbb{C}} C(\mathbb{T}))$ by \[ \pi_1\left( \begin{bmatrix} a & b \\
c & d
\end{bmatrix}\right) = \begin{bmatrix} a & bz_1 \\ c \overline{z_1} & d
\end{bmatrix} \] where $z_1$ is the usual generator for
$C(\mathbb{T})$ in the first copy of $C(\mathbb{T}) \subseteq
C(\mathbb{T}) \free{\mathbb{C}} C(\mathbb{T})$. A routine
calculation shows that $\pi_1$ is a $*$-representation.

Next define $\pi_2: M_2 \rightarrow M_2(C(\mathbb{T})
\free{\mathbb{C}} C(\mathbb{T}))$ by \[ \pi_2\left( \begin{bmatrix} a & b \\
c & d
\end{bmatrix}\right) = \begin{bmatrix} \frac{a + d - c\overline{z_2} + b
z_2}{2} & \frac{a-d-c\overline{z_2}+bz_2}{2} \\
\frac{a-d+c\overline{z_2} - bz_2}{2} & \frac{
a+d+c\overline{z_2}+bz_2}{2} \end{bmatrix} \] where $z_2$ is the
usual generator for $C(\mathbb{T})$ in the second copy of
$C(\mathbb{T}) \subseteq C(\mathbb{T}) \free{\mathbb{C}}
C(\mathbb{T})$. Again, a routine calculation shows that $\pi_2$ is a
$*$-representation.

Now $\pi_1\left(\begin{bmatrix} a & 0 \\ 0 & a \end{bmatrix}\right)
= \begin{bmatrix} a & 0 \\ 0 & a \end{bmatrix} =  \pi_2\left(
\begin{bmatrix} a & 0 \\ 0 & a \end{bmatrix}\right)$ and hence there is a
$*$-representation $\pi_1 \free{} \pi_2: M_2 \free{\mathbb{C}} M_2
\rightarrow M_2( C(\mathbb{T}) \free{\mathbb{C}} C(\mathbb{T}))$.
Notice that \[ \begin{bmatrix} z_1 & 0 \\ 0 & 0 \end{bmatrix} =
\pi_1\left(\begin{bmatrix} 0 & 1 \\ 0 & 0 \end{bmatrix}\right)
\pi_2\left( \begin{bmatrix} 1 & 0 \\ 0 & -1 \end{bmatrix}\right)
 \] and \[
\begin{bmatrix} z_2 & 0 \\ 0 & 0 \end{bmatrix} =
\pi_1\left(\begin{bmatrix} 1 & 0 \\ 0 & 0 \end{bmatrix}\right)
\pi_2\left(
\begin{bmatrix} 0 & 2 \\ 0 & 0 \end{bmatrix}\right) \pi_1\left( \begin{bmatrix}
1 & 0 \\ 0 & 0 \end{bmatrix}\right) \] and hence the non-exact
subalgebra $
\begin{bmatrix} C(\mathbb{T}) \free{\mathbb{C}} C(\mathbb{T}) & 0 \\
0 & 0 \end{bmatrix}$ is contained as a subalgebra in the image of
$\pi_1 \free{}\pi_2$.  It follows that $ M_2 \free{\mathbb{C}} M_2$
can not be exact. \end{proof}

It is not hard to see that, in the previous proof, the mapping
$\pi_1 \free{} \pi_2$ is not one-to-one.  This follows since $
\mathbb{C}^2 \free{\mathbb{C}} \mathbb{C}^2$ is a subalgebra of $
M_2 \free{\mathbb{C}} M_2$, but the image of $ \mathbb{C}^2
\free{\mathbb{C}} \mathbb{C}^2$ under the mapping $\pi_1 \free{}
\pi_2$ is finite dimensional.  However it is well known, see
\cite[Example IV.1.4.2]{Blackadar} that $ \mathbb{C}^2
\free{\mathbb{C}} \mathbb{C}^2$ is isomorphic to
\[ \left\{ \begin{bmatrix} f_{1,1}(t) & f_{1,2}(t) \\ f_{2,1}(t) &
f_{2,2}(t) \end{bmatrix}: f_{i,j} \in C([0,1]), f_{1,2}(0) =
f_{2,1}(0) = f_{1,2}(1) = f_{2,1}(1) = 0 \right\}. \]

\section{Algebras of the form $M_j \free{D} M_k \free{D} M_l$}

We first note that $M_j \free{D} M_k \free{D} M_l = M_j \free{D} M_l
\free{D} M_k$ and hence if any two of $j,k$, or $l$ give rise to
amalgamation diagrams with minimum value greater than or equal to
$3$, then the algebra $M_j \free{D} M_k \free{D} M_l$ is not exact.

\begin{theorem} The algebra $M_2 \free{\mathbb{C}^2} M_2 \free{
\mathbb{C}^2} M_2$ is isomorphic to $M_2(C(\mathbb{T})
\free{\mathbb{C}} C(\mathbb{T}))$ and hence is not exact.
\end{theorem}

\begin{proof}  By Example \ref{m2on}, the algebra $M_2
\free{\mathbb{C}^2} M_2$ is isomorphic to $M_2 \otimes
C(\mathbb{T})$.  Further there is a canonical $*$-isomorphism \[
\pi: (M_2 \free{\mathbb{C}^2} M_2) \free{M_2} (M_2
\free{\mathbb{C}^2} M_2) \rightarrow M_2 \free{\mathbb{C}^2} M_2
\free{\mathbb{C}^2} M_2.\]  Now assume that $\pi_i: C(\mathbb{T})
\otimes M_2 \rightarrow A$ are unital $*$-representations satisfying
$ \pi_1 (1 \otimes d) = \pi_2(1 \otimes d) $ for all $ d \in M_2$.
For $ a \in C(\mathbb{T})$ we define $ \sigma_i (a) = \pi_i(a
\otimes 1)$. Then there exists $\sigma_1 \free{} \sigma_2 :
C(\mathbb{T}) \free{\mathbb{C}} C(\mathbb{T}) \rightarrow A$.
Further, if we set $\sigma: M_2 \rightarrow A$ by $\sigma(d) =
\sigma_1(1 \otimes d)$ then we know that $ \sigma(d) \sigma_1(a) =
\sigma_1(a) \sigma(d)$ and $\sigma(d) \sigma_2(a) = \sigma_2(a)
\sigma(d)$ for all $a \in C(\mathbb{T})$ and $ d \in M_2$ and hence
$ \sigma_1 \free{}\sigma_2 (x) \sigma(d) = \sigma(d) \sigma_1
\free{}\sigma_2 (x)$ for all $x \in C(\mathbb{T}) \free{\mathbb{C}}
C(\mathbb{T})$ and $d \in D$. It follows by the universal property
of the tensor product that there exists $\tau: C(\mathbb{T})
\free{\mathbb{C}} C(\mathbb{T}) \otimes M_2 \rightarrow A$ extending
the canonical inclusions of $C(\mathbb{T}) \otimes M_2$ into
$(C(\mathbb{T}) \free{\mathbb{C}} C(\mathbb{T})) \otimes M_2$. Hence
$(M_2 \free{\mathbb{C}^2} M_2) \free{M_2} (M_2 \free{\mathbb{C}^2}
M_2)$ is isomorphic to $(C(\mathbb{T}) \free{\mathbb{C}}
C(\mathbb{T})) \otimes M_2$ which is not exact since it contains a
copy of $C(\mathbb{T}) \free{\mathbb{C}} C(\mathbb{T})$. \end{proof}

\begin{corollary} Let $D$ be a diagonal subalgebra of $M_2$, then $M_2 \free{D} M_2
\free{D} M_2$ is not exact. \end{corollary}

\begin{proposition} Let $D$ be a unital diagonal subalgebra of $M_j,
M_k$ and $M_l$.  The algebra $M_j \free{D} M_k \free{D} M_l$ is not
exact. \end{proposition}

\begin{proof}  If the dimension of $D$ is greater than or equal to
$3$ then by Proposition \ref{mjmk} the algebra can not be exact, so
we look only at the case that $\dim D \leq 2$.

If $\dim D = 2$, then again we can assume without loss of generality
that in the threefold amalgamation diagram for the free product that
in any given row at most one box is not equal to $1$.  Thus at least
one of $j,k,$ or $l$ must equal $2$.  So without loss of generality
assume that we have $l = 2$ and we are in the case of $M_j
\free{\mathbb{C}^2} M_k \free{\mathbb{C}^2} M_2$.  Now, as in
Proposition \ref{submjmk} we can see that there is a copy of $M_2
\free{\mathbb{C}^2} M_2$ inside $M_j \free{\mathbb{C}^2} M_k$, and
applying \cite{Armstrong-Dykema-Exel-Li:2003} again we have that $
M_2 \free{\mathbb{C}^2} M_2 \free{\mathbb{C}^2} M_2$ is a subalgebra
of $M_j \free{\mathbb{C}^2} M_k \free{\mathbb{C}^2}  M_l$ and hence
the latter is not exact.

The case of $\dim D = 1$ now follows by Corollary \ref{subalgebra}.
\end{proof}

\section{Free products with no amalgamation and some nonunital amalgamations}

We know by applying Proposition \ref{mjmk} and Proposition
\ref{subalgebra} that the following is true.

\begin{proposition} The algebra $M_j \free{} M_k$ is not exact for
any $k,j \geq 2$. \end{proposition}

We now focus on the case in which the diagonal subalgebra $D$
contains the identity of $M_j$ but not that of $M_k$.  In this case
the amalgamation diagram is
of the form \begin{align*} & M_j: \boxed{j_1} \boxed{j_2} \cdots \boxed{j_m} & \\
& M_k:\boxed{k_1} \boxed{k_2} \cdots \boxed{k_m} \boxed{0} \\
\end{align*} where $m$ is the dimension of $D$.  We will write
$k(D)$ for the integer given by $ k - \sum_{i=1}^m k_i$

\begin{theorem} Let $D$ be a  unital diagonal subalgebra of $M_j$,
where $D$ is a diagonal subalgebra of $M_k$ which does not contain
the unit of $M_k$.  Then the algebra $M_j \free{D} M_k$ is exact if
and only if $M_j \free{D} M_{k-k(D)}$ is exact in which case $M_j
\free{D} M_k$ is nuclear.
\end{theorem}

\begin{proof} Clearly, since $M_j \free{D} M_{k-k(D)}$ is a
subalgebra of $M_j \free{D} M_k$ if the former is not exact neither
is the latter.  We will focus on the case in which $M_j \free{D}
M_{k-K(D)}$ is exact.  This breaks down into two cases.

Case 1 ($\dim D = 1$): In this case, either $j = 1$ which is
trivial, or $k-k(D) = 1$ which puts us in the context of Example
\ref{mntensora}.

Case 2 ($\dim D = 2$): In this case the subalgebra $M_j \free{D}
M_{k-K(D)}$ is a directed graph algebra, see Proposition
\ref{graph}.  The corresponding directed graph has two vertices $\{
v_1, v_2\}$ and $ j-1$ edges from $v_1$ to $v_2$ and $ k-k(D)-1$
edges from $v_2$ to $v_1$.  Create a new graph $G$ by adding a
vertex $v_3$ and $k-k(D)$ edges from $v_2$ to $v_3$.  We claim that
the algebra $C^*(G)$ is isomorphic to $M_j \free{D} M_k$ and hence
the algebra is nuclear.

Let $\{ E, P \}$ be the Cuntz-Krieger system given by the generators
for the graph $C^*$-algebra $M_j \free{D} M_{k-k(D)}$.  Now look at
the associated Cuntz-Krieger system \[ \left\{ E \cup \{e_{j,k}: 1
\leq j \leq k-1\}, P \cup \left\{ \sum_{m= k-k(D)+1}^{k}e_{m,m}
\right\} \right\},\] where $e_{i,j} \in M_{k} \subset M_j \free{D}
M_{k}$. Notice that this new Cuntz-Krieger system generates $M_{j}
\free{D} M_{k}$ as a $C^*$-algebra and hence a standard result for
graph algebras, \cite[Proposition 1.21]{Raeburn:2006}, gives an onto
representation $: \pi: C^*(G) \rightarrow M_j \free{D} M_{k}$.  Now
since graph algebras are nuclear the algebra $M_j \free{D} M_k$ is
nuclear.\end{proof}

Finally we can make some progress on the general case. We know that
there is a canonical onto $*$-representation $\pi: M_3 \free{D} M_3
\rightarrow M_3 \free{\mathbb{C}^3} M_3$ for any diagonal subalgebra
$D$ and hence $M_3 \free{D} M_3$ is not exact for any diagonal
subalgebra $D$.  We have also seen that $M_2 \free{} M_2$ and the
free product with unital amalgamation, $M_2 \free{\mathbb{C}} M_2$,
are not exact.

Now write the amalgamation diagram for $M_j \free{D} M_k$ as
\begin{align*} & M_j:\boxed{j_1} \boxed{j_2}  \cdots \boxed{j_m} \boxed{0} \\
& M_k:\boxed{k_1} \boxed{k_2}  \cdots \boxed{k_m} \boxed{0}.
\end{align*}

\begin{proposition} Let $D$ be a non-unital diagonal subalgebra of $M_j$ and
$M_k$. If $\dim D \geq 2$ then $M_j \free{D} M_k$ is not exact. If
either $ j-\sum j_i$ or $k-\sum k_i$ is greater than or equal to $2$
then $M_j \free{D} M_k$ is not exact.  \end{proposition}

\begin{proof}  We deal first with $\dim D \geq 2$.  Notice that
there will be an embedding of $M_3$ into $M_j$ and $M_k$ so that the
subalgebra will have amalgamation diagram \begin{align*} & M_3:
\boxed{1} \boxed{1} \boxed{0} \\ & M_3:\boxed{1} \boxed{1} \boxed{0}
\end{align*} which will have as a quotient the non-exact algebra
$M_3 \free{\mathbb{C}^3} M_3$ and hence $M_j \free{D} M_k$ will not
be exact.

For the other situation we notice that there will be a subalgebra of
the form $\mathbb{C} \free{} (\mathbb{C} \oplus \mathbb{C})$.  This
non-unital $C^*$-algebra satisfies \[ (\mathbb{C} \free{}
(\mathbb{C} \oplus \mathbb{C}))^1 \cong \left(\mathbb{C} \oplus
\mathbb{C}\right) \free{\mathbb{C}} \left( \mathbb{C} \oplus
\mathbb{C} \oplus \mathbb{C}\right).\]  The latter algebra is
isomorphic to $C^*(\mathbb{Z}_2) \free{\mathbb{C}} C^*(\mathbb{Z}_3)
\cong C^*(\mathbb{Z}_2 \free{} \mathbb{Z}_3)$ which contains a copy
of $C^*(\mathbb{Z}\free{}\mathbb{Z})$ which is not exact.  It
follows that since the unitization of $\mathbb{C} \free{}
(\mathbb{C} \oplus \mathbb{C})$ is not exact the algebra is not
either and hence $M_j \free{D} M_k$ is not exact.
\end{proof}

The only case that remains is the free product $M_2
\free{\mathbb{C}} M_k$ with amalgamation diagram \begin{align*} &
M_2:\boxed{1} \boxed{0} \\ & M_k:\boxed{k-1} \boxed{0}.
\end{align*}  We do not, as yet, have a satisfactory answer for this
situation.

\bibliographystyle{plain}

\end{document}